\documentclass[12pt]{amsart}
\usepackage[russian,english]{babel}
\usepackage{amssymb}
\usepackage{tikz}
\usepackage{tikz-cd}
\usepackage{tkz-euclide}

\newtheorem{theorem}{Theorem}

\newtheorem{proposition}{Proposition}

\newtheorem{remark}{Remark}

\newcommand{\Z}{{\mathbb Z}}

\newcommand{\R}{{\mathbb R}}
\newcommand{\CN}{{\mathbb C}}

\title{Complements of  discriminants of real boundary singularities}

\begin{document}

\author{M.A.~Gudiev}\address{Higher School of Economics, Moscow, Russia} \email{ magudiev@gmail.com}

\begin{abstract}
We study the topology of the complements of discriminants of simple real boundary singularities by counting the connected components of these sets and assigning to them certain topological characteristics. Results of this paper serve as a generalization of those recently acquired by Vassiliev \cite{VVA} for ordinary function singularities.
\end{abstract}


\maketitle

\section{Introduction}

In \cite{ArnADE} Arnol'd found a deep connection between simple function singularities and the Dynkin diagrams $A_k$, $D_k$ and $E_k$ and their corresponding Weyl groups. Later in \cite{ArnBCF} he expanded this classification to the diagrams of types $B_k$, $C_k$ and $F_4$, which turned out to be closely connected to simple singularities of functions on manifolds with boundary. The present work expands upon the work of V.A.Vassiliev (see \cite{VVA}) by considering simple boundary singularities of types $B_k$, $C_k$ and $F_4$.

For a function singularity $f:(\R^n, 0) \to (\R, 0)$, $df|_0 = 0$ and its arbitrary smooth deformation $F:(\R^n \times \R^k, 0) \to (\R,0)$, which can be viewed as a family of functions $f_\lambda = F(-,\lambda):\R^n \to \R$, s.t. $f_0 = f$, we can define the set $\Sigma = \Sigma(F) \subset \R^k$, called the \textit{(real) discriminant} of $F$ to be the set of all parameters $\lambda$, s.t. $f_\lambda$ has a zero critical value. In cases we will consider this will be an algebraic subset of codimension one in $\R^k$, which divides a small neighborhood of the origin in $\R^k$ into several parts. The following theorem holds for standard versal deformations of real simple function singularities.
\begin{theorem}[E.Looijenga, \cite{Lo}]
All connected components of the complements of the real discriminant varieties of standard versal deformations of simple real function singularities are contractible.
\end{theorem} 
This theorem implies that the topology of the set $\R^k \setminus \Sigma$ is completely defined by the number of connected components of this set. We will prove a similar result for singularities $B_k$  and $C_k$. The author believes that this should also be true in the case $F_4$, but this is yet to be proven.

The topology and combinatorics of such complements was described by V.D. Sedykh in his works \cite{Sed,SedBook} for simple singularities with Milnor number $\leq 6$ (see \cite{Sed}, Theorems 2.8 and 2.9 for the numbers of local components for singularities $D_4$, $D_5$, $D_6$ and $E_6$). Recently (see \cite{VVA}) Vassiliev fully described the topology of the complement $\R^k \setminus \Sigma$ for simple real function singularities and their versal deformations by listing the number of local components in each case and assigning a certain topological characteristic to each of them. We will enumerate the local components of the complements for $B_k$, $C_k$ and $F_4$ boundary singularities and assign to them similar topological invariants.

Any simple function singularity, up to stable equivalence, can be realized as a function $f:\R^2 \to \R$ in two variables and has a versal deformation of dimension $\mu$, where $\mu$ is the \textit{Milnor number} of $f$. For any parameter $\lambda \in \R^\mu$ of a versal deformation of a simple singularity, consider the set of lower values $W(\lambda) = \{x \in \R^2 | f_\lambda(x) \leq 0\}$. These sets can go to infinity along several \textit{asymptotic sectors} the number of which stays the same for any $\lambda$ for a given versal deformation. We say that two sets of lower values $W(\lambda_1)$ and $W(\lambda_2)$ are \textit{topologically equivalent} if there exists an orientation-preserving homeomorphism of $\R^2$ which sends $W(\lambda_1)$ to $W(\lambda_2)$ but doesn't permute the asymptotic sectors. Naturally, if two parameters $\lambda_1$ and $\lambda_2$ lie in the same connected component of $\R^\mu \setminus \Sigma$ then the corresponding sets of lower values are topologically equivalent. The main theorem of \cite{VVA} states the converse is also true:

\begin{theorem}[Vassiliev,\cite{VVA}]
If $\lambda_1$ and $\lambda_2$ are non-discriminant points of the parameter space $\R^\mu$ of a versal deformation of a simple function singularity, and the corresponding sets  $W(\lambda_1), W(\lambda_2)$ are topologically equivalent, then $\lambda_1$ and $\lambda_2$ belong to the same component of $\R^\mu \setminus \Sigma$.
\end{theorem}

Our main goal will be to prove a similar statement for simple boundary singularities, and describe the connected components of the complement of the discriminant and their corresponding sets of lower values. 

\section{Notions and Definitions}

Assuming the reader is familiar with basic notions of singularity theory (for a classic reference see \cite{AVG1,AVG2}) we will only describe analogues of the usual constructions for the case of boundary function singularities. Another good reference is \cite{Lyash}.

To avoid ambiguity, further in the text we will use the term ordinary singularity for singularities of functions on manifolds without boundary.

Consider the space $\R^n$ with a fixed hyperplane $\{x = (x_1, \ldots, x_n) \in \R^n|x_1 = 0\}$, which will act as a germ of an $n$-dimensional manifold with boundary.

A \textit{(real) boundary function singularity} is a germ of a function $f: (\R^n,0) \to (\R,0)$ on a manifold with boundary such that $0$ is a critical point of the restriction of $f$ onto the boundary, i.e.
$$
\frac{\partial f}{\partial x_i} \bigg|_{x=0} = 0, \; i =2,\ldots,n
$$ 

We'll call two boundary singularities $f_i$, $i = 1,2$ equivalent if there exists a local diffeomorphism $\varphi: \R^n \to \R^n$ preserving the boundary s.t. $\varphi^* f_2 = f_1$. 

Hence the classification problem can be formulated in terms of describing the orbits of action of the group $Loc_{B}(\R^n)$ of local diffeomorphisms preserving the boundary on the space of function germs. The \textit{modality} of a singularity $f$ is the minimal number $m$ such that a sufficiently small neighborhood of $f$ in the space of germs can be covered by a finite number of no more than $m$-parametric orbits of the action of $Loc_B(\R^n)$. A singularity of modality $0$ will be called \textit{simple}. As was mentioned earlier, simple boundary singularities are classified by diagrams $B_k$, $C_k$ and $F_4$.

A \textit{gradient ideal} of an ordinary singularity is the ideal generated by its partial derivatives:

$$
I_f = (\frac{\partial f}{\partial x_1}, \ldots, \frac{\partial f}{\partial x_n} )
$$

which in case of boundary singularities is defined as

$$
I_{f|x_1} = (x_1 \frac{\partial f}{\partial x_1}, \ldots, \frac{\partial f}{\partial x_n}).
$$

\textit{Local algebra} of a germ $f$ is defined as the factor algebra

$$
Q_{f} = \R[[x_1,\ldots,x_n]]/I_f
$$
and its boundary analog as

$$
Q_{f|x_1} = \R[[x_1,\ldots,x_n]]/I_{f|x_1}.
$$

It's easy to see that for equivalent singularities the corresponding local algebras are isomorphic, so this construction gives a powerful invariant of singularities.

The \textit{Milnor number} of an ordinary singularity is defined as $\mu = \dim Q_f$, and in case of a boundary singularity as $\mu = \dim Q_{f|x_1}$. A classic result is that a germ $f$ has finite Milnor number iff $f$ is an isolated singularity. For boundary singularities we also define two additional numbers

$$
\mu_0 = \dim \R[[x_1,\ldots,x_n]]/(\partial f/\partial x_2 |_{x_1 = 0}, \ldots, \partial f/\partial x_n |_{x_1 = 0})
$$

$$
\mu_1 = \dim \R[[x_1,\ldots,x_n]]/(\frac{\partial f}{\partial x_1}, \ldots, \frac{\partial f}{\partial x_n}).
$$
The number $\mu_0$ is the Milnor number of $f$ viewed as a function on the boundary $\{x_1 = 0\}$ and $\mu_1$ is the Milnor number of $f$ viewed as an ordinary germ in $\R^n$. It's easy to see that $\mu = \mu_1 + \mu_0$. Singularities for which $\mu_1 = 0$ will be called purely boundary.

All the above definitions also work if we take $\CN$ as the base field and consider holomorphic functions instead of smooth ones.

The classification of boundary singularities also includes ordinary ones: from a function $f:\R^n \to \R$ we can obtain a purely boundary singularity $\tilde{f}(x_0,x) = x_0 + f(x)$, for a manifold $\R^{n+1}$ with boundary $\{(x_0,x) | x_0 = 0\}$. Thus the simple boundary function singularities also include ones of types $A_k$, $D_k$ and $E_k$. The remaining types, specific to the boundary case, are listed in table 1.

\begin{table}
\begin{center}
\begin{tabular}{c c c}
Notation & Normal form & Number of components\\ 
$B^+_{2k}, k \geq 1$ & $x^{2k} + y^2$ & $(k+1)^2$\\
$B^-_{2k}, k \geq 1$ & $x^{2k} - y^2$ & $(k+1)^2$\\
$\pm B_{2k + 1}, k \geq 1$ & $\pm(x^{2k +1} + y^2)$ & $(k+1)(k+2)$\\
$C^+_{2k}, k \geq 1$ & $xy + y^{2k}$ & $(k+1)^2$\\
$C^-_{2k}, k \geq 1$ & $xy - y^{2k}$ & $(k+1)^2$\\
$\pm C_{2k + 1}, k \geq 1$ & $\pm (xy + y^{2k + 1})$ & $(k+1)(k+2)$\\
$F_4$ & $\pm x^2 + y^3$ & 8
 \end{tabular} \end{center}
\caption{Normal forms of real simple boundary singularities in two variables $(x,y)$, with boundary given by $x = 0$.}
\label{t1}
\end{table}

A \textit{deformation} of a germ $f(x)$ is a function
$$
F(x,\lambda):(\R^n \times \R^l, 0) \to (\R,0)
$$
s.t. $F(x,0) = f(x)$. The space $\R^l$ is called the base of deformation $F$. We call a deformation $F$ \textit{versal} if any other deformation can be induced from it, meaning for any $G:\R^n\times \R^k \to \R$, $G(x,0) = f(x)$ we can find 
\begin{enumerate}
\item a germ of a smooth map $\psi:\R^k \to \R^l$ 
\item and a diffeomorphism $\eta(x,\xi)$ of $\R^n$, $\xi\in \R^k$, which depends smoothly on $\xi$ and s.t. $\eta(x,0) = x$, so that
\end{enumerate}
$$
G(x, \xi) = F(\eta(x,\xi),\psi(\xi))
$$

Geometrically a deformation is just a germ of a surface at the point $f$ in the space of germs, and it is versal iff this germ is transversal to the orbit of $f$ under the action of $Loc(\R^n)$. The deformation is called \textit{miniversal}, if the dimension of the parameter space is minimal. For a boundary function singularity $f$ with local algebra $Q_{f|x_1}$ and (finite) Milnor number $\mu$ the (mini)versal deformation can be obtained by setting 
$$
F(x,\lambda) = f(x) + \sum_{i = 1}^\mu \lambda_i f_i(x)
$$
where $f_i$ form a basis of the local algebra $Q_{f|x_1}$. Of course, the same construction holds in the ordinary case. Hence, the miniversal deformations for simple boundary singularities can be given as follows:
\begin{eqnarray}
B_\mu & \qquad & f_\lambda(x,y) = x^\mu \pm y^2 + \lambda_1 x^{\mu-1} + \ldots + \lambda_\mu \label{bvd} \\
C_\mu & \qquad & f_\lambda(x,y) = xy \pm y^\mu + \lambda_1 y^{\mu-1} + \ldots + \lambda_\mu \label{cvd} \\
 F_4 & \qquad & f_\lambda(x,y) = \pm x^2 + y^3 + \lambda_1 x + \lambda_2 y + \lambda_3 xy + \lambda_4 \label{f4vd}
\end{eqnarray}

Given a (mini)versal deformation $F(x,\lambda) = f_\lambda(x)$ of a boundary singularity we define the \textit{discriminant variety} as the subset $\Sigma \subset \R^\mu$ of the parameter space consisting of all singular parameter values. A parameter $\lambda \in \R^\mu$ is called \textit{non-singular}, if 

\begin{enumerate}
\item $0$ is a regular value of the function $f_\lambda$
\item the zero level set $V(\lambda) = f^{-1}_\lambda(x)$ is transversal to the boundary.
\end{enumerate}

In case of ordinary singularities only the first condition is required. For a boundary singularity the discriminant consists of two parts corresponding to these conditions. We will denote by $\Sigma_0$ the component corresponding to the first condition, and by $\Sigma_1$ the one corresponding to the second. The discriminant of an ordinary singularity is an irreducible affine variety, and it follows that in the boundary case the components $\Sigma_i$ are also irreducible. 

In her note \cite{Shcher} I. G. Shcherbak introduced the notion of \textit{decomposition} of a boundary singularity $f$. The decomposition is defined as a pair (type of $f$ as an ordinary singularity, type of the restriction $f|_{\{x_1=0\}})$. Naturally, any boundary singularity possesses such a decomposition, moreover there exists an involution on the set of boundary singularities that swaps these types. The singularities that go into each other under this involution are called \textit{dual}. In particular, the singularity $B_\mu$ has a decomposition $(A_{\mu-1}, A_1)$, $C_\mu$ the decomposition $(A_1, A_{\mu-1})$ and $F_4$ the decomposition $(A_2, A_2)$. This means that the singularities $B_\mu$ and $C_\mu$ are dual (and $B_2 = C_2$ is self-dual), and $F_4$ is self-dual.

For a boundary singularity $f$ that decomposes into $(f_0, f_1)$ the sets $\Sigma_0$ and $\Sigma_1$ are diffeomorphic to the discriminants $\Sigma_{f_0}$ and $\Sigma_{f_1}$ multiplied by euclidean spaces of suitable dimensions. This will be useful later in the study of the discriminant set of $F_4$.

As it was mentioned earlier, the following statement holds:

\begin{proposition}
All connected components of the complements of the real discriminant varieties of versal deformations 1-2 of singularities $B_\mu$, $C_\mu$ are contractible.
\end{proposition}

Recall that by $W(\lambda) = f_\lambda^{-1}((-\infty, 0])$ we denote the set of lower values for $\lambda$.
Given a deformation of a boundary singularity, we say that two sets of lower values $W(\lambda_1)$ and $W(\lambda_2)$ are topologically equivalent if there exists an orientation preserving homeomorphism of $\R^n$ sending one set to the other, which preserves the boundary and does not permute the asymptotic sectors of $f$ and its restriction to the boundary $f|_{x_1 = 0}$. The number of asymptotic sectors is equal to $0$ in case $B^+_{2k}$, $1$ in cases $B_{2k+1}$ and $F_4$ and $2$ in cases $C_\mu$ and $B^-_{2k}$. Now we are ready to formulate the main theorem:

\begin{theorem}
The numbers of components of the complements of the real discriminant varieties of deformations 1-3 are listed in Table 1, and each of these components is uniquely defined by the topological type of corresponding set of lower values.
\end{theorem}

The proofs of theorems 3 and 4 in cases $B_\mu$ and $C_\mu$ will be presented in the next section, and the proof in case $F_4$ in sections 4 and 5.

\section{Cases $B_\mu$ and $C_\mu$}

As it was mentioned earlier, the discriminant $\Sigma$ consists of two components $\Sigma_0$ and $\Sigma_1$ each defined by an appropriate condition. For $B_\mu$ singularities the first condition means that the polynomial $x^\mu + \lambda_1 x^{\mu-1} + \ldots + \lambda_\mu $ can only have simple roots, and the second condition means that $0$ is not a root of this polynomial. For $C_\mu$ these conditions interchange.

Denote by 
$$
h_\lambda(x) = \pm x^\mu + \lambda_1 x^{\mu-1} + \ldots + \lambda_\mu,
$$
then the deformation can be expressed as $f_\lambda(x,y) = h_\lambda(x) \pm y^2$, and the equation for the zero set is given by $y^2 = -h_\lambda(x)$. Hence the zero set is symmetric with respect to the line $\{y = 0\}$ and consists of some number of ovals, each intersecting the line $\{y=0\}$ at neighboring pairs of roots of $h_\lambda$, and no more than two non-compact components.

Any connected component of the complement $\R^\mu \setminus \Sigma$ is completely defined by the configuration of roots of the polynomial $h_\lambda(x)$. If $\lambda$ is a non-discriminant parameter, then, since the polynomial $h_\lambda$ must only have simple non-zero roots, it can only have $p$ negative and $q$ positive roots, and since they all should be simple the parity of $p+q$ should coincide with that of $\mu$. For any $\lambda$ such that $h_\lambda$ has $p$ negative and $q$ positive roots there is an obvious path to any other $\lambda'$ with equal numbers of negative and positive roots. Moreover, we can construct a homotopy contracting the whole connected component of $\lambda$ to any of its points. It's also easy to see that the numbers $p, q$ uniquely define the topological type of $W(\lambda)$.

In case $C_\mu$ the zero set is given by the equation
$$
xy = -h_\lambda(y),
$$
hence, if $\lambda$ is non-discriminant, it consists of the graph of the function $x = -h_\lambda(y)/y$ (notice that by changing the equation this way we don't lose the solution $y=0$, as otherwise $\lambda$ would be a discriminant parameter). This function has a vertical asymptote $y=0$ and, as before, the topological type of $W(\lambda)$ and the connected component of $\lambda$ is completely defined by the numbers $p$ and $q$ of negative and positive roots of $h_\lambda$. 

The above description gives proofs of proposition 1 and theorem 3 in cases $B_\mu$ and $C_\mu$.

\section{The case $F_4$}

In this section we will obtain $6$ out of possible $8$ types of sets of lower values for $F_4$ and study a certain hyperplane section of the discriminant $\Sigma$.

For brevity we will denote the versal deformation of $F_4$ by
$$
f_\lambda(x,y) = x^2 + y^3 + a x + b y + c xy + d.
$$
The non-transversality condition for points of $\Sigma_1$ means the polynomial $y^3 + b y + d$ has a non-simple root. A direct calculation gives the following equation: 
$$
27d^2 + 4b^3 = 0,
$$
meaning $\Sigma_1$ is a direct product of a cusp in the plane $(b, d)$ and a plane $\R^2$ spanned by coordinates $a, c$. 

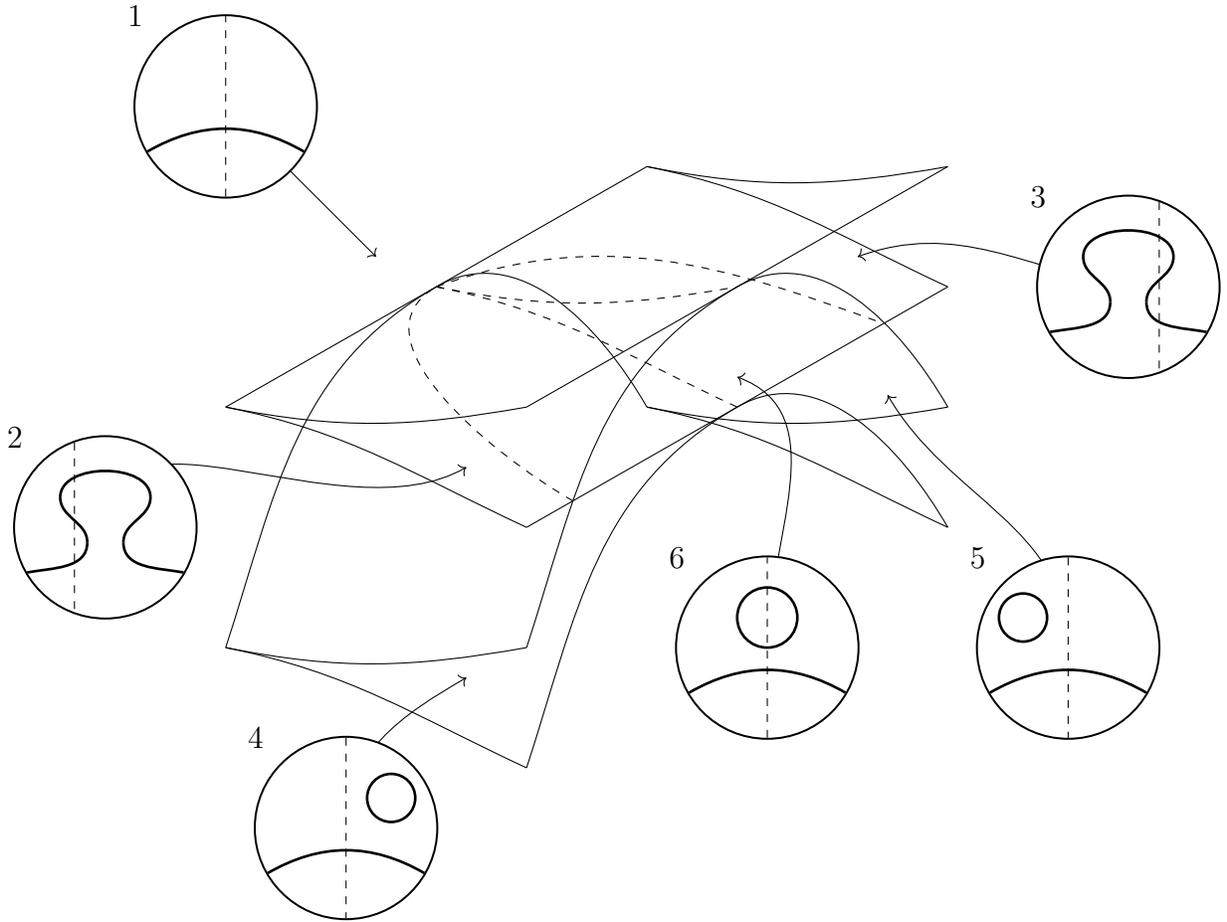
\begin{figure}
\begin{tikzpicture}[scale=0.8]
\draw (0,0) -- (7,4);
\draw (0,0) to[out=-11.3,in=190] (5,0);
\draw (0,0) to[out=-11.3, in=155] (5,-2);
\draw (5,0) -- (12,4);
\draw (7,4) to[out=-11.3,in=190] (12,4);
\draw (7,4) to[out=-11.3, in=155] (12,2);
\draw (5,-2) -- (12,2);
\draw (0,-4) to[out=-11.3,in=190] (5,-4);
\draw (0,-4) to[out=-11.3, in=155] (5,-6);
\draw (7,0) to[out=-11.3,in=190] (12,0);
\draw (7,0) to[out=-11.3, in=155] (12,-2);
\draw (0,-4) to[out=72, in=210] (3.5,2);
\draw (3.5,2) to[out=30, in=120] (7,0);
\draw (5,-4) to[out=72, in=210] (8.5,2);
\draw (8.5,2) to[out=30, in=120] (12,0);
\draw (5,-6) to[out=72, in=210] (8.5,0);
\draw (8.5,0) to[out=30, in=120] (12,-2);
\draw [dashed](3.5,2) to[out=-11.3,in=190] (8.5,2);
\draw [dashed](3.5,2) to[out=-11.3, in=155] (8.5,0);
\draw [dashed](5.77,-1.56) to[out=150,in=210] (3.5,2);
\draw [dashed](3.5,2) to[out=20, in=160] (10.93,1.39);
\draw [line width = 1.5pt] (0,5) circle (1.5);
\draw (-1.5,6.5) node {1};
\draw [->] (0,5) -- (2.5,2.5);
\fill [white] (0,5) circle (1.5);
\draw [dashed] (0,3.5) -- (0,6.5);
\draw [line width = 1pt] (-1.3,4.25) to[out=30,in=150] (1.3,4.25);
\draw [line width = 1.5pt] (9,-4) circle (1.5);
\draw (9 - 1.5, -4 + 1.5) node {6};
\draw [->] (9,-4) to[out = 90, in = -20] (8.5,0.5);
\fill [white] (9,-4) circle (1.5);
\draw [dashed] (9,-5.5) -- (9,-2.5);
\draw [line width = 1pt] (9,-3.5) circle (0.5);
\draw [line width = 1pt] (7.7,-4.75) to[out=30,in=150] (10.3,-4.75);
\draw [line width = 1.5pt] (-2,-2) circle (1.5);
\draw (-2-1.5, -2 + 1.5) node {2};
\draw [->] (-2,-2) to[out = 90, in = 210] (4,-1);
\fill [white] (-2,-2) circle (1.5);
\draw [dashed] ({-2 + 1.5*cos(110)},{-2 + 1.5*sin(110)}) -- ({-2 + 1.5*cos(250)},{-2 + 1.5*sin(250)});
\draw [line width = 1pt] (-0.7,-2.75) to[out=170,in=270] (-1.7,-2.25);
\draw [line width = 1pt] (-3.3,-2.75) to[out=10,in=270] (-2.3,-2.25);
\draw [line width = 1pt] (-1.7,-2.25) to[out=90,in=270] (-1.25,-1.5);
\draw [line width = 1pt] (-2.3,-2.25) to[out=90,in=270] (-2.75,-1.5);
\draw [line width = 1pt] (-1.25,-1.5) to[out=90, in=90] (-2.75,-1.5);
\draw [line width = 1.5pt] (15,2) circle (1.5);
\draw (15 - 1.5, 2 + 1.5) node {3};
\draw [->] (15,2) to[out = 170, in = 20] (10.5,2.5);
\fill [white] (15,2) circle (1.5);
\draw [dashed] ({15 + 1.5*cos(70)},{2 + 1.5*sin(70)}) -- ({15 + 1.5*cos(290)},{2 + 1.5*sin(290)});
\draw [line width = 1pt] (-0.7 + 17,-2.75 + 4) to[out=170,in=270] (-1.7 + 17,-2.25 + 4);
\draw [line width = 1pt] (-3.3 + 17,-2.75 + 4) to[out=10,in=270] (-2.3 + 17,-2.25 + 4);
\draw [line width = 1pt] (-1.7 + 17,-2.25 + 4) to[out=90,in=270] (-1.25 + 17,-1.5 + 4);
\draw [line width = 1pt] (-2.3 + 17,-2.25 + 4) to[out=90,in=270] (-2.75 + 17,-1.5 + 4);
\draw [line width = 1pt] (-1.25 + 17,-1.5 + 4) to[out=90, in=90] (-2.75 + 17,-1.5 + 4);
\draw [line width = 1.5pt] (2,-7) circle (1.5);
\draw (2 - 1.5, -7 + 1.5) node {4};
\draw [->] (2,-7) to[out = 90, in = 210] (4,-4.5);
\fill [white] (2,-7) circle (1.5);
\draw [dashed] (2,-5.5) -- (2,-8.5);
\draw [line width = 1pt] (2.75,-6.5) circle (0.4);
\draw [line width = 1pt] (0.7,-7.75) to[out=30,in=150] (3.3,-7.75);
\draw [line width = 1.5pt] (14,-4) circle (1.5);
\draw (14 - 1.5, -4 + 1.5) node {5};
\draw [->] (14,-4) to[out = 90, in = 300] (11,0.2);
\fill [white] (14,-4) circle (1.5);
\draw [dashed] (14,-5.5 + 3) -- (14,-8.5 + 3);
\draw [line width = 1pt] (13.25,-3.5) circle (0.4);
\draw [line width = 1pt] (0.7 + 12,-7.75 + 3) to[out=30,in=150] (3.3 + 12,-7.75 + 3);
\end{tikzpicture}
\caption{The $3$-dimensional section of the discriminant variety of $F_4$ and corresponding zero sets}
\label{F4ds}
\end{figure}

As we will see later, the equation for $\Sigma_0$ is also computable, however it turns out to be quite complex. However, if we restrict our deformation to dimension $3$ by setting $c = 0$, we will get a nice section of $\Sigma_0$ by the plane $\{c = 0\}$, in which the equation will have the following form:
$$
27 \left(d + \frac{a^2}{4}\right)^2 + 4b^3 = 0
$$
meaning $\Sigma'_0 = \Sigma_0 \cap \{c = 0\}$ is a cuspidal edge bent along the parabola given by equations $d = -a^2/4, b = 0$. As depicted in fig. 1 edges of $\Sigma'_0$ and $\Sigma'_1$ are tangent at the origin, and these sets are tangent along the cusp lying in the plane $\{a = 0\}$. It's easy to see that the local components of the set $\R^3 \setminus \Sigma' = \R^3 \setminus (\Sigma'_0 \cup \Sigma'_1)$ (here $\R^3$ denotes the hyperplane $\{c=0\}$ of the parameter space) are all contractible, and the number of these components is equal to $6$. As also depicted in fig. 1, corresponding sets of lower values are all distinct and correspond to different topological types of the real elliptic curve with respect to the boundary
$$
x^2 + y^3 + a x + b y + d = 0.
$$

Concrete realizations of these sets are easy to obtain by taking values of the parameter $\lambda = (a,b,0,d)$ to lie in the corresponding component.

\begin{figure}
\begin{tikzpicture}[scale=1]
\draw (5,3) -- (7,3) -- (7,1) -- (5,1) -- (5,3);
\draw[dashed] (5.75,3) -- (5.75,1);
\draw (5,2) to[out=30,in=90] (6.25,2) to[out=270,in=90] (5.25,1.5) to[out=270,in=180] (7,1.5);
\draw (6.65, 2.3) circle (0.25);
\draw (6,0.5) node {7};
\draw (7.5,3) -- (9.5,3) -- (9.5,1) -- (7.5,1) -- (7.5,3);
\draw[dashed] (8.75,3) -- (8.75,1);
\draw (9.5,2) to[out=150,in=90] (8.25,2) to[out=270,in=90] (9.25,1.5) to[out=270,in=0] (7.5,1.5);
\draw (7.85, 2.3) circle (0.25);
\draw (8.5,0.5) node {8};
\draw (0,0) -- (2,0) -- (2,-2) -- (0,-2) -- (0,0);
\draw[dashed] (1,0) -- (1,-2);
\draw (0,-1) to[out=30,in=150] (2,-1);
\draw (0.5, -1.5) circle (0.3);
\draw (1,-2.5) node {9};
\draw (0 + 2.5,0) -- (2 + 2.5,0) -- (2 + 2.5,-2) -- (0 + 2.5,-2) -- (0 + 2.5,0);
\draw[dashed] (1 + 2.5,0) -- (1 + 2.5,-2);
\draw (0 + 2.5,-1) to[out=30,in=150] (2+ 2.5,-1);
\draw (3.5, -1.5) circle (0.3);
\draw (1 + 2.5,-2.5) node {10};
\draw (0 + 5,0) -- (2 + 5,0) -- (2 + 5,-2) -- (0 + 5,-2) -- (0 + 5,0);
\draw[dashed] (1 + 5.25,0) -- (1 + 5.25,-2);
\draw (5,-0.5) to[out=0,in=90] (6.5,-0.5) to[out=270,in=90] (5.25,-1.25) to[out=270,in=180] (7,-1.75);
\draw (5.75, -1.3) circle (0.25);
\draw (1 + 5,-2.5) node {11};
\draw (0 + 7.5,0) -- (2 + 7.5,0) -- (2 + 7.5,-2) -- (0 + 7.5,-2) -- (0 + 7.5,0);
\draw[dashed] (8.25,0) -- (8.25,-2);
\draw (7.5,-0.25) to[out=0,in=90] (9.25,-0.5) to[out=270,in=90] (7.75,-1.5) to[out=270,in=180] (9.5,-1.75);
\draw (8.75, -0.5) circle (0.25);
\draw (8.5,-2.5) node {12};
\draw (0 + 10,0) -- (2 + 10,0) -- (2 + 10,-2) -- (0 + 10,-2) -- (0 + 10,0);
\draw[dashed] (11.25,0) -- (11.25,-2);
\draw (10,-0.25) to[out=0,in=90] (11.75,-0.5) to[out=270,in=90] (10.25,-1.5) to[out=270,in=180] (12,-1.75);
\draw (10.5, -0.615) circle (0.25);
\draw (11,-2.5) node {13};
\draw (0 + 12.5,0) -- (2 + 12.5,0) -- (2 + 12.5,-2) -- (0 + 12.5,-2) -- (0 + 12.5,0);
\draw[dashed] (13.25,0) -- (13.25,-2);
\draw (12.5,-0.25) to[out=0,in=90] (14.25,-0.5) to[out=270,in=90] (12.75,-1) to[out=270,in=180] (14.5,-1);
\draw (14, -1.5) circle (0.25);
\draw (13.5,-2.5) node {14};
\end{tikzpicture}
\caption{Possible remaining topological types of sets of lower values for $F_4$. Notice that for each set in the bottom row, except for №10, the one reflected through the boundary is also a possible type.}
\label{F4ds}
\end{figure}
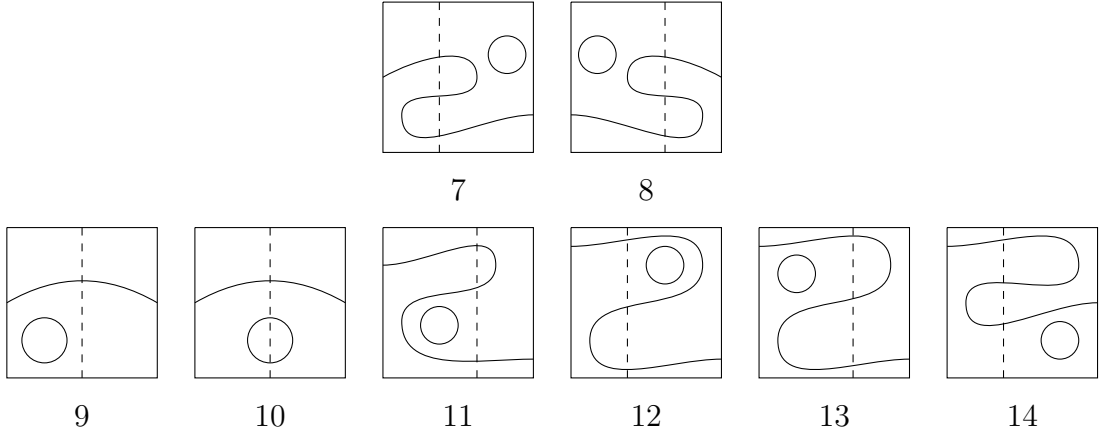

\section{Further calculations for $F_4$}

In order to complete the proof of the main theorem we need to describe the topology of the set $\R^4 \setminus \Sigma$ for $F_4$. In this section we will find out which remaining topological types of sets of lower values can be realized through the versal deformation and calculate the homology of the complement.

\subsection{Remaining topological types of $W(\lambda)$  for $F_4$}

For the versal deformation $f_\lambda(x,y)$ the corresponding zero set $f_\lambda^{-1}(0)$ is either a non-compact line going to infinity along the $x$-axis, or a union of such a line with a compact oval. Notice also, that since the polynomial $f_\lambda(0,y)$ in $y$ has degree $3$ the number of points of intersection of the zero set with the boundary $\{x=0\}$ is equal to $1$ or $3$ if $\lambda$ is a non-discriminant parameter. Thus the possible remaining topological types of sets of lower values look like ones listed in fig.2.

However, only the sets $7$ and $8$ can be realized as ones coming from the versal deformation of $F_4$, and as we'll see later that these are the only ones remaining. 

Concrete realizations of the types $7$ and $8$ can be obtained as follows. First, a direct computation shows that the component of the set $\Sigma_0 \cap \Sigma_1$ (which will be further denoted by $\Xi_0$) corresponding to functions $f_\lambda$ that have a Morse critical point with zero critical value which lies on the boundary, can be parametrized through $c$ and $d$ the following way:

$$
f_\lambda(x,y) = x^2 + y^3 - c \sqrt[3]{\frac{d}{2}} x - 3\sqrt[3]{\frac{d}{4}} y + cxy + d.
$$

We first take a sufficiently small $d > 0$ and $c \neq 0$ (the sign of $c$ dictates what type of set of lower values, $7$ or $8$, will be obtained), so that the root of multiplicity $2$ of the polynomial $f_\lambda(0,y)$ in $y$ would be greater than the remaining root. This produces a zero set shown in the leftmost picture of fig.3. Notice that a substitution of the form $x \mapsto x - \varepsilon$ can be realized through our deformation by an appropriate change of parameters $d$ and $a$, hence we can move the boundary away from the crossing to obtain the middle picture of fig.3. Finally, by subtracting a small constant $\delta > 0$ from our function we can remove the crossing so that the curve splits into two separate components, hence we obtain the sets $7$ and $8$.

\begin{figure}
\begin{tikzpicture}[scale=2]
\draw[line width=1pt](0,0) -- (2,0) -- (2,2) -- (0,2) -- (0,0);
\draw[dashed,line width=1pt] (1,0) -- (1,2);
\draw[line width=0.9pt](0,0.5) to[out=20,in=220] (1,1) to[out=40,in=300] (1.75, 1.5) to[out=120,in=60] (1,1) to[out=240,in=90] (0.75,0.5) to[out=270,in=220] (2,0.5);
\draw[->] (2.2,1) -- (2.8,1);
\draw[line width=1pt] (3,0) -- (5,0) -- (5,2) -- (3,2) -- (3,0);
\draw[dashed, line width=1pt] (3.87,0) -- (3.87,2);
\draw[line width=0.9pt] (3,0.5) to[out=20,in=220] (4,1);
\draw[line width=0.9pt]  (4,1) to[out=40,in=300] (4.75, 1.5);
\draw[line width=0.9pt] (4.75, 1.5) to[out=120,in=60] (4,1);
\draw[line width=0.9pt] (4,1) to[out=240,in=90] (3.75,0.5);
\draw[line width=0.9pt] (3.75,0.5) to[out=270,in=220] (5,0.5);
\draw[->] (5.2,1) -- (5.8,1);
\draw[line width=1pt] (6,0) -- (8,0) -- (8,2) -- (6,2) -- (6,0);
\draw[dashed, line width=1pt] (6.87,0) -- (6.87,2);
\draw[line width=0.9pt] (6,0.5) to[out=10,in=90] (7,1);
\draw[line width=0.9pt]  (7,1) to[out=270,in=90] (6.5, 0.5);
\draw[line width=0.9pt] (6.5, 0.5) to[out=270,in=210] (8,0.5);
\draw[line width=0.9pt] (7.5,1.4) circle (0.25);
\end{tikzpicture}
\caption{Construction of the remaining sets of lower values}
\label{78}
\end{figure}
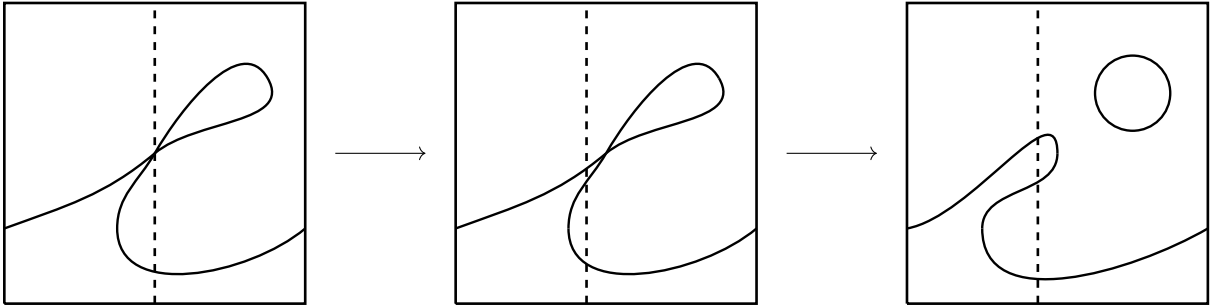

\subsection{Homological calculations for the complement}

One way to calculate the number of connected components of the set $\R^4 \setminus \Sigma$ is to study its reduced cohomology groups $\tilde{H}^i(\R^4 \setminus \Sigma)$. By Alexander duality ($\tilde{H}^*$ denotes the reduced cohomology group, $\bar{H}_*$ - the Borel-Moore homology) we have
$$
\tilde{H}^i(\R^4 \setminus \Sigma) \simeq \bar{H}_{3-i}(\Sigma)
$$

We will prove the following theorem, which will imply the completeness of the lists of topological types given in this and previous section.

\begin{theorem}
The (reduced) homology group $\bar{H}_i(\Sigma;\Z_2)$ is isomorphic to $(\Z_2)^{7}$ if $i=3$ and is trivial otherwise.
\end{theorem}

\begin{remark} 
Such homology groups can be studied by standard methods developed by Vassiliev (see \cite{V88, Vbook}), which were used in \cite{VVA} in order to compute cohomology of the complements of the discriminant varieties of $D_\mu$ singularities. We will use a different approach.
\end{remark}

\begin{proof}
The conditions on the parameters of $\Sigma_0$ yield the following system of polynomial equations:

$$
\begin{cases}
f_\lambda(x,y) = x^2 + y^3 + a x + b y + c xy + d = 0\\
\frac{\partial f_\lambda}{\partial x} = 2x + a + cy = 0 \\
\frac{\partial f_\lambda}{\partial y} = 3y^2 + cx + b = 0
\end{cases}
$$

To get an equation on $\lambda$ we can substitute $x$ by a linear function in $y$ using the second equation. We obtain a system of two polynomial equations, for which we can then write down the resultant in $y$ to eliminate it. Notice that taking the resultant does not add any imaginary solutions for this system, as this would imply that the deformation $f_\lambda$, as a function of complex variable, has two distinct conjugate critical points with critical value $0$, which is impossible for the singularity $A_2$, which is the type of $F_4$ as an ordinary singularity.

As it was mentioned earlier, the equation for $\Sigma_1$ has the form 
$$
27d^2 + 4b^3 = 0,
$$
so we obtain a system of two polynomial equations for the intersection $\Sigma_0 \cap \Sigma_1$. The solution of this system consists of two $2$-dimensional irreducible components $\Xi_0$ and $\Xi_1$, each homeomorphic to $\R^2$. The component $\Xi_0$ is comprised of parameters $\lambda$ corresponding to deformations $f_\lambda$ which have an ordinary Morse critical point with critical value $0$ lying on the boundary and the closure of such points, and the component $\Xi_1$ corresponds to deformations for which the zero set is non-transversal to the boundary and which have an ordinary Morse critical point with critical value $0$ outside the boundary (and their closure).

Once again, we can calculate the intersection $\Xi_0 \cap \Xi_1$. Turns out, $\Xi_0$ intersects $\Xi_1$ along two curves $\Psi_0$ and $\Psi_1$ which intersect transversely at the origin and are both homeomorphic to $\R^1$. The curve $\Psi_0$ corresponds to deformations which have a Morse critical point at the origin, such that one of the branches of the curve $f^{-1}_\lambda(0)$ is tangent to the boundary. The curve $\Psi_1$ is comprised of parameters $\lambda$ for which the deformation has a critical point of type $A_2$ lying on the boundary.

Recall that the two components of the discriminant $\Sigma$ of a boundary singularity can be obtained from the discriminants of the ordinary singularities into which it decomposes. As $F_4$ has the decomposition $(A_2,A_2)$, the components $\Sigma_0$ and $\Sigma_1$ are both diffeomorphic to a cusp multiplied by $\R^2$, meaning $\Sigma_i$ are both homeomorphic to $\R^3$. Now we can calculate the Borel-Moore homology of $\Sigma$ by first applying the Mayer-Vietoris long exact sequence to the decomposition $\Sigma_0 \cap \Sigma_1 = \Xi_0 \cup \Xi_1$ and after that to $\Sigma = \Sigma_0 \cup \Sigma_1$ (the $\Z_2$ coefficients are omitted):

\begin{tikzcd}
\ldots \rar & \bar{H}_j(\Psi_0 \cup \Psi_1) \rar & \bar{H}_j(\Psi_0) \oplus \bar{H}_j(\Psi_1) \rar & \bar{H}_j(\Xi_0 \cup \Xi_1) \rar & \ldots 
\end{tikzcd}

Since the one-point compactification of $\Psi_0 \cup \Psi_1$ is homotopy equivalent to a bouquet of $3$ circles, we get that $\bar{H}_1(\Psi_0 \cup \Psi_1) \simeq (\Z_2)^3$ and $0$ in other dimensions. Hence, from the long exact sequence we get

$$
\bar{H}_j(\Xi_0 \cup \Xi_1) = \begin{cases}
(\Z_2)^5 & \text{for} \; j = 2, \\
0 & \text{otherwise.}
\end{cases}
$$

Now applying this calculation to the exact sequence

\begin{tikzcd}
\ldots \rar & \bar{H}_j(\Sigma_0 \cap \Sigma_1) \rar & \bar{H}_j(\Sigma_0) \oplus \bar{H}_j(\Sigma_1) \rar & \bar{H}_j(\Sigma) \rar & \ldots 
\end{tikzcd}

we get that 

$$
\bar{H}_j(\Sigma) = \begin{cases}
(\Z_2)^7 & \text{for} \; j = 3, \\
0 & \text{otherwise,}
\end{cases}
$$
which completes the proof.
\end{proof}

\end{document}